\documentclass[reqno, 11pt]{amsart}
\usepackage{amsfonts,latexsym}
\usepackage{amsmath}
\usepackage{amscd}
\usepackage{float,amsmath,amssymb,mathrsfs,bm,multirow,graphics}
\usepackage[dvips]{graphicx}
\usepackage[percent]{overpic}
\usepackage[all]{xy}

\newcommand{\nequation}{\setcounter{equation}{0}}
\renewcommand{\theequation}{\mbox{\arabic{equation}}}
\newcommand{\R}{{\Bbb R}}

\newcommand{\proofbegin}{\noindent{\it Proof.\,\,}}
\newcommand{\proofend}{\hfill$\Box$\bigskip}
\newcommand{\llangle}{\langle\!\langle}
\newcommand{\rrangle}{\rangle\!\rangle}

\newcommand{\Diff}{\text{\upshape Diff}}

\newcommand{\Diffmu}{\Diff_{\mu}}
\newcommand{\DiffM}{\Diff(M)}

\newcommand{\DiffmuM}{\Diffmu(M)}

\newcommand{\id}{\text{\upshape id}}

\newcommand{\ad}{\text{\upshape ad}}
\newcommand{\Ad}{\text{\upshape Ad}}
\newcommand{\Laplacian}{\Delta}

\newcommand{\met}{}

\DeclareMathOperator{\diver}{div}

\DeclareMathOperator{\im}{im}
\DeclareMathOperator{\grad}{grad}

\allowdisplaybreaks

\newtheorem{theorem}{Theorem}
\newtheorem{proposition}[theorem]{Proposition}

\newtheorem{lemma}[theorem]{Lemma}

\newtheorem{remark}[theorem]{Remark}

\newtheorem{figuretext}{Figure}

\input epsf
\title[A two-component geodesic equation]{\sc A two-component geodesic equation on a space of constant positive curvature}
\author{Jonatan Lenells}
\email{Jonatan\_Lenells@baylor.edu}

\author{Zhao Yang}
\email{Zhao\_Yang@baylor.edu}

\address{Department of Mathematics, Baylor University, One Bear Place \#97328, Waco, TX 76798, USA.}

\begin{document}

\begin{abstract}
\noindent
We propose a new two-component geodesic equation with the unusual property that the underlying space has constant positive curvature. In the special case of one space dimension, the equation reduces to the two-component Hunter-Saxton equation.
\end{abstract}

\maketitle




\section{Introduction}\nequation
Many basic equations in physics, including the Euler equations of motion of a rigid body, and the Euler equations of fluid dynamics of an inviscid incompressible fluid, can be viewed as geodesic equations on a Lie group endowed with a one-sided invariant metric \cite{A1966}. In the context of rigid bodies, the group is $SO(3)$; in the context of incompressible fluids, it is the group $\Diff_\mu(\R^3)$ of volume-preserving diffeomorphisms. Other geodesic equations include Burgers' equation ($\Diff(S^1)$ with a right-invariant $L^2$ metric), the KdV equation (Virasoro group with a right-invariant $L^2$ metric), and the Camassa-Holm equation ($\Diff(S^1)$ with a right-invariant $H^1$ metric), see \cite{V2008} for a survey.

In the recent work \cite{KMLP2011}, the following equation was introduced and studied:
\begin{align}\label{h1doteq}
d\diver{u}_t = - d\left(\iota_u d\diver(u) + \frac{1}{2}\diver(u)^2\right),
\end{align}
where $u$ is a time-dependent vector field on a compact Riemannian manifold $M$, $d$ denotes the exterior derivative, and $\iota$ denotes the interior product.
Two of the remarkable properties of (\ref{h1doteq}) are: (a) It is the geodesic equation on a space of constant positive curvature. Indeed, equation (\ref{h1doteq}) is the equation for geodesic flow on the space of right cosets $\Diff(M)/\Diff_\mu(M)$ endowed with the right-invariant metric given at the identity by
\begin{align*}
\llangle u,v\rrangle = \frac{1}{4} \int_M \diver(u) \diver(v) d\mu,
\end{align*}
and this space has constant curvature equal to $1/\mu(M) > 0$.
(b) It is an integrable evolution equation in any number of space dimensions. Indeed, it is shown in \cite{KMLP2011} that (\ref{h1doteq}) admits a complete set of independent integrals
in involution.

In this paper we propose the following two-component generalization of (\ref{h1doteq}):
\begin{align}\label{2h1doteq}
\begin{cases}
d \diver{u}_t = - d\left(\iota_u d\diver(u) + \frac{1}{2}\diver(u)^2 - \frac{1}{2}\rho^2\right),
	\\
 \rho_t = -\diver(\rho u),
 \end{cases}
\end{align}	
where $\rho$ is a time-dependent real-valued function on $M$.
We will show that (\ref{2h1doteq}) is the geodesic equation on the homogeneous space $[\DiffM \circledS C^\infty(M)]/\DiffmuM$ endowed with the right-invariant metric given at the identity by
\begin{align} \label{2h1dotmetric}
\llangle (u_1, u_2), (v_1, v_2) \rrangle = \frac{1}{4} \int_M [\diver(u_1)\diver(v_1) + u_2v_2] d\mu,
\end{align}
and that this space has constant positive curvature.

Our search for a two-component generalization of (\ref{h1doteq}) was motivated by the special case when $M$ is the unit circle, i.e. $M = S^1$. In this case, equation (\ref{h1doteq}) reduces to the Hunter-Saxton equation \cite{HS1991}
\begin{align}\label{hs}
u_{txx} = -\Bigr(uu_{xx} + \frac{1}{2}u_x^2\Bigr)_x,
\end{align}
which admits the following integrable two-component generalization, see \cite{CI2008, W2009}:
\begin{align}\label{2hs}
 \begin{cases}
   u_{txx} = -\left(uu_{xx} + \frac{1}{2}u_x^2 - \frac{1}{2}\rho^2\right)_x,
  	\\
  \rho_t = -(\rho u)_x,	
  \end{cases}
\end{align}
Equation (\ref{hs}) describes the geodesic flow on a sphere of constant positive curvature \cite{KM2003, Lsphere} and it was recently observed that the two-component version (\ref{2hs}) also is the geodesic equation on a space of constant positive curvature \cite{K2011}.

Since the three special cases (\ref{h1doteq}), (\ref{hs}), and (\ref{2hs}) of the system (\ref{2h1doteq}) are all integrable, we expect (\ref{2h1doteq}) to also be integrable, but the proof of this remains open.

Our main result is stated in section \ref{mainresultsec} and its proof is presented in section \ref{proofsec}. 

\section{Main result}\label{mainresultsec}
Let $M$ be a compact connected $n$-dimensional Riemannian manifold and let $\DiffM$ denote the space of orientation-preserving diffeomorphisms of $M$. We let $G$ denote the semidirect product $G = \DiffM \circledS C^\infty(M)$ with multiplication given by
$$(\varphi, f)(\psi, g) = (\varphi\circ\psi, g + f \circ \psi), \qquad (\varphi, f), (\psi, g) \in G.$$
Moreover, we let $H = \DiffmuM$ denote the subgroup of $\DiffM$ of volume-preserving diffeomorphisms. The group $H$ acts on $G$ by
$$\psi \cdot (\varphi, f) = (\psi \circ \varphi, f), \qquad \psi \in H, \; (\varphi, f) \in G,$$
and our main interest lies in the homogeneous space $G/H$ of right cosets.

We use right-invariance to extend (\ref{2h1dotmetric}) to a (degenerate) metric $\llangle \cdot, \cdot \rrangle$ on $G$. Thus,
$$\llangle U, V \rrangle_{(\varphi, f)} = \frac{1}{4} \int_M \bigl[\diver(U_1 \circ \varphi^{-1})\diver(V_1 \circ \varphi^{-1}) + (U_2 \circ \varphi^{-1})(V_2 \circ \varphi^{-1})\bigr] d\mu, $$
where $U = (U_1, U_2)$ and Ê$V = (V_1, V_2)$ are elements of $T_{(\varphi, f)}G$ and $d\mu$ denotes the volume form on $M$ induced by the metric.
It is straightforward to verify that this metric descends to a Riemannian metric on $G/H$ (see section \ref{proofsec} for details), and we can now state our main result.

\begin{theorem}\label{mainth}
The sectional curvature of $G/H$ equipped with the metric $\llangle \cdot, \cdot \rrangle$ is constant and equal to $1/\mu(M) > 0$. Moreover, the corresponding Euler equation for the geodesic flow is the two-component equation (\ref{2h1doteq}).
\end{theorem}

\begin{remark}\upshape
As far as regularity assumptions and function spaces are concerned, the geometry of equation (\ref{2h1doteq}) can be developed in a number of different settings, including the $H^s$ and $C^n$ settings (see \cite{EKL2011} for further details in a similar situation). Since these issues are of no consequence for the considerations here, we have chosen for simplicity to use the notation of the smooth category.
\end{remark}

\section{Proof}\label{proofsec}
We will use the following notations:
\begin{itemize}
\item $e = (\id, 0)$ will denote the identity element of $G$.
\item $\mathfrak{X}(M)$ will denote the space of vector fields on $M$.
\item For $p \geq 0$, $\Omega^p(M)$ will denote the space of $p$-forms on $M$.
\item $\langle \cdot, \cdot \rangle$ will denote the metric on $M$. $\langle \cdot, \cdot \rangle$ will also denote the induced inner product on $p$-forms, as well as the natural pairing between elements of $TM$ and $T^*M$.
\item If $X,Y$ are vector fields on $M$, $[X,Y]$ will denote their Lie bracket, given locally by $[X,Y] = DY \cdot X - DX \cdot Y$.
\item If $u = (u_1, u_2)$ and $v = (v_1, v_2)$ are elements of $T_eG \simeq \mathfrak{X}(M) \times C^\infty(M)$, $[u,v] \in T_eG$ will denote the following commutator:
$$[u,v] = \begin{pmatrix} [u_1, v_1] \\ dv_2(u_1) - du_2(v_1) \end{pmatrix}.$$
\end{itemize}

\subsection{The metric on $G/H$}
We will first show that the metric $\llangle \cdot, \cdot \rrangle$ descends to $G/H$.
According to the general theory of homogeneous spaces (see Proposition 3.1 in Chapter X of \cite{KN1969}), $\llangle \cdot, \cdot \rrangle$ decends to $G/H$ provided that
\begin{equation} \label{descent}
\llangle \ad_w u, v\rrangle + \llangle u, \ad_w v\rrangle = 0 \quad \text{for all} \quad u,v\in T_eG \; \text{and}\; w\in T_eH,
\end{equation}
where we view $H \simeq \DiffmuM \times \{0\} \subset G$ as a subgroup ofÊ $G$.
Using that the adjoint action in $T_eG$ is given by
$$\ad_{u}v
= -[u,v], \qquad u, v \in T_eG,$$
it is straightforward to verify (\ref{descent}): Given elements $u=(u_1, u_2),v=(v_1, v_2)$, and $w = (w_1, w_2)$, we have
\begin{align}\label{checkdescent}
\llangle \ad_w u, v\rrangle & + \llangle u, \ad_w v\rrangle
= \frac{1}{4}\int_M \big( \diver{[u_1, w_1]} \diver{v_1} +  \diver{u_1} \diver{[v_1,w_1]}   \big) d\mu
	\\ \nonumber
&+ \frac{1}{4}\int_M \big( (dw_2(u_1) - du_{2}(w_1))v_2 + u_2 (dw_2(v_1) - dv_2(w_1))  \big) d\mu.
\end{align}
If $w \in T_e H$ then $\diver{w_1} = 0$ and $w_2 = 0$, and the right-hand side of (\ref{checkdescent}) equals
\begin{align*}
&\frac{1}{4}\int_M  \Big[\big( \langle u_1 , \grad \diver{w_1} \rangle -  \langle w_1, \grad \diver{u_1} \rangle \big) \diver{v_1}
	\\
&\hspace{1.3cm} + \diver{u_1} \big( \langle v_1, \grad\diver{w_1} \rangle - \langle w_1, \grad\diver{v_1} \rangle \big)  \Big] d\mu
- \frac{1}{4}\int_M \diver(u_2v_2w_1) d\mu
	\\
= & -\frac{1}{4}\int_M \diver(w_1  \diver{v_1} \diver{u_1})d\mu  = 0.
\end{align*}
This shows that $\llangle \cdot, \cdot \rrangle$ descends to $G/H$.

\subsection{The Euler equation}
We next show that the Euler equation associated with the metric $\llangle \cdot, \cdot \rrangle$ on $G/H$ is the two-component equation (\ref{2h1doteq}). It follows from the general theory developed by Arnold \cite{A1966} that the Euler equation is given by
\begin{align}\label{Eulereq}
  u_t = B(u, u),
\end{align}
where $u(t)$ is a curve in $T_eG$ and the bilinear operator $B:T_eG \times T_eG \to T_eG$ is defined by the condition
\begin{align}\label{Bcondition}
  \llangle B(u, v), w \rrangle = \llangle u, [v, w] \rrangle, \qquad u,v,w \in T_eG.
\end{align}
In order to determine $B$, we define the operator $A:\mathfrak{X}(M) \to \Omega^1(M)$ by $Av = d\delta v^\flat$, where $\delta$ denotes the codifferential on $M$ and $\flat:TM \to T^*M$ denotes the musical isomorphism with inverse $\sharp:T^*M \to TM$ induced by the metric $\langle \cdot, \cdot \rangle$. Recalling that $\diver(X) = -\delta X^\flat$ for a vector field $X$, the condition (\ref{Bcondition}) can be written as
\begin{align}\label{Bcondition2}
\int_M \langle AB_1(u,v), w_1 \rangle d\mu + \int_M B_2(u,v) w_2 d\mu
= \int_M \langle Au_1, [v,w]_1 \rangle d\mu
+ \int_M u_2 [v,w]_2d\mu,
\end{align}
where $B_1$ and $B_2$ denote the two components of $B$. Using the general formula
\begin{align}\label{lieformula}
  [v_1,w_1]^\flat =  \diver(w_1) v_1^\flat - \diver(v_1) w_1^\flat -\delta (v_1^\flat \wedge w_1^\flat),
\end{align}
we see that the right-hand side of (\ref{Bcondition2}) equals
\begin{align*}
 & \int_M \langle Au_1, \diver(w_1) v_1^\flat - \diver(v_1) w_1^\flat  - \delta(v_1^\flat \wedge w_1^\flat)  \rangle d\mu
  + \int_M u_2(dw_2(v_1) - dv_2(w_1)) d\mu
  	\\
&= \int_M \langle  - d\langle Au_1, v_1 \rangle - \diver(v_1) Au_1 - \iota_{v_1}dAu_1, w_1 \rangle d\mu	
- \int_M \bigl(w_2\diver(u_2v_1) + \langle u_2dv_2, w_1 \rangle\bigr) d\mu. 	
\end{align*}
Thus, since $w$ is arbitrary and $dAu_1 = 0$, we find
\begin{align}\nonumber
  AB_1(u,v) & = -d\langle Au_1, v_1 \rangle - \diver(v_1) Au_1 - u_2 dv_2,
  	\\ \label{Bexpression}
  B_2(u,v) & = -\diver(u_2 v_1).	
\end{align}
Substituting this expression for $B = (B_1, B_2)$ into (\ref{Eulereq}), we find the two-component system (\ref{2h1doteq}).

\subsection{The operator $A$}
Equation (\ref{Bexpression}) only fixes the value of $B_1$ up to an element of $\ker{A}$. This is a reflection of the fact that the inner product (\ref{2h1dotmetric}) is degenerate along $T_eH$ in $T_eG$. Accordingly, the solution of the Euler equation (\ref{Eulereq}) is not uniquely determined in $T_eG$ by the initial data, but descends to a well-defined curve in $T_e(G/H)$. 

The following lemma, which will be needed for the computation of the curvature, clarifies some properties of $A$ and $A^{-1}$. In particular, it shows that the kernel of $A$ equals $T_eH$ and consists of all divergence-free vector fields on $M$. 

\begin{lemma}\label{Alemma}
The operator $A:\mathfrak{X}(M) \to \Omega^1(M)$ defined by $Av = d\delta v^\flat = - d\diver{v}$ satisfies
\begin{align}\label{imageA}
   \im{A} = \left\{df \, | \, f \in C^\infty(M) \right\}
\end{align}   
and
\begin{align}\label{kerA}
  \ker{A} = \{v \in \mathfrak{X}(M)\, | \, \diver{v} = 0\}.
\end{align}  
In particular, the map $\diver{A^{-1}d}:C^\infty(M) \to C^\infty(M)$ is well-defined and is given by
\begin{align}\label{deltaAinvd}
  \diver{A^{-1}df} = -f + \frac{1}{\mu(M)}\int_M f d\mu.
\end{align}  
\end{lemma}
\proofbegin
The inclusion $\subset$ in (\ref{imageA}) is obvious. 
On the other hand, if $f \in C^\infty(M)$, then standard properties of the Laplacian $\Laplacian = d\delta + \delta d$ on a compact manifold (see \cite{R1997}) imply that $\ker \Laplacian$ consists of the constant functions on $M$ and that there exists $g \in C^\infty(M)$ such that $\Laplacian g = f - \frac{1}{\mu(M)}\int_M f d\mu \in (\ker\Laplacian)^\perp$. Letting $v = (dg)^\sharp$, we find $-\diver{v} = f - \frac{1}{\mu(M)}\int_M f d\mu$ and so
$$Av = -d\diver{v} = df.$$
This proves (\ref{imageA}).
In order to prove (\ref{kerA}), we first note that if $\diver{v} = 0$, then $Av = 0$. Coversely, suppose $Av = 0$. Then
$$0 = \int_M \langle d\delta v^\flat, v^\flat \rangle d\mu 
= \int_M \langle \delta v^\flat, \delta v^\flat \rangle d\mu$$ 
and we find $\delta v^\flat = -\diver{v} = 0$. This proves (\ref{kerA}).

It is clear from (\ref{imageA}) and (\ref{kerA}) that the map $\diver{A^{-1}d}:C^\infty(M) \to C^\infty(M)$ is well-defined. 
Letting $h = \diver{A^{-1}df}$, we find
$$dh = -AA^{-1}df = -df.$$
Thus, $h = -f + c$, where $c$ is a constant. Noting that 
$$0 = \int_M h d\mu = -\int_M f d\mu + c \mu(M),$$
we infer that
$$c = \frac{1}{\mu(M)}\int_M f d\mu.$$
This proves (\ref{deltaAinvd}). 
\proofend

\begin{figure}
\begin{center}
  \begin{overpic}[width=.6\textwidth]{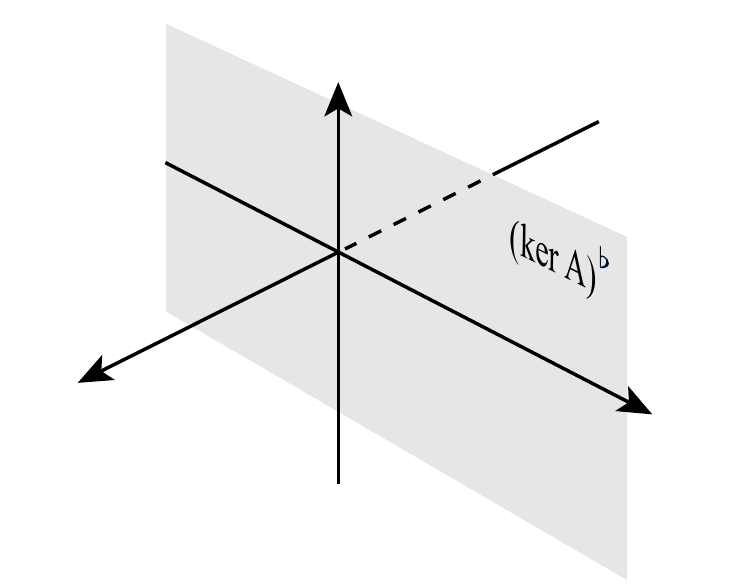}
    \put(42,73){ $\mathcal{H}^1$}
    \put(-12,24){ $d\Omega^0(M)$}
     \put(90,20){ $\delta \Omega^2(M)$}
 \end{overpic}
\begin{figuretext} \label{hodge.pdf}
The Hodge decomposition of $T_e^*\Diff(M)$ and $\ker{A}$.
\end{figuretext}
\end{center}
\end{figure}

\begin{remark}\upshape
In terms of the Hodge decomposition
\begin{equation}\label{hodgedecomposition}
T_e^*\Diff(M) = d\Omega^0(M) \oplus \delta \Omega^2(M) \oplus \mathcal{H}^1,
\end{equation}
where $\mathcal{H}^k$ denotes the space of harmonic $k$-forms on $M$, we have $\ker{A} = T_eH = (\delta \Omega^2(M) \oplus \mathcal{H}^1)^\sharp$, see Figure \ref{hodge.pdf}. It follows that the value of $B_1$ can be fixed by requiring that $B_1(u,v) \in \{(df)^\sharp | f \in C^\infty(M)\}$. However, we will not need to do this. 
\end{remark}

\subsection{Constant positive curvature}
It remains to prove that $(G/H, \llangle \cdot, \cdot \rrangle)$ has constant sectional curvature equal to $1/\mu(M)$. Letting $R$ denote the curvature tensor on $G/H$, this is equivalent to proving that
\begin{align}\label{Rexpression}
  \llangle R(u,v)v, u \rrangle =  \frac{\llangle u, u \rrangle \llangle v, v \rrangle - \llangle u, v\rrangle^2}{\mu(M)}, \qquad u,v \in T_eG.
\end{align}
We will use the following curvature formula which is derived in the appendix:
\begin{align}\label{curvatureformula}
\llangle R(u,v)v, u\rrangle
 = \llangle \delta, \delta \rrangle + \llangle [u,v], \beta \rrangle - \frac{3}{4}\llangle [u,v], [u,v]\rrangle - \llangle B(u,u), B(v,v) \rrangle,
\end{align}
where $\delta$ and $\beta$ are defined by
$$\delta = \frac{1}{2}(B(u,v) + B(v,u)), \qquad \beta = \frac{1}{2}(B(u,v) - B(v,u)).$$

\begin{remark}\upshape
1. Formula (\ref{curvatureformula}) coincides formally with the curvature formula for a Lie group equipped with a left-invariant metric given in Appendix 2 of \cite{Abook}.
In the appendix, we show that this formula remains valid in the case of a right-invariant metric on a homogeneous space.

2. Although $B_1$ is only defined up to an element of $\ker{A}$, Lemma \ref{Alemma} implies that the right-hand side of (\ref{curvatureformula}) is well-defined. Indeed, we will see that $A^{-1}$ only enters (\ref{curvatureformula}) via the combination $\diver{A^{-1}df}$ and by (\ref{deltaAinvd}) this operator is uniquely defined. 
\end{remark}

Our proof of (\ref{Rexpression}) will be a long computation using (\ref{curvatureformula}) together with the expression (\ref{Bexpression}) for $B$. First note that
\begin{align}\nonumber
\delta & = -\frac{1}{2}\begin{pmatrix} A^{-1}[d\langle Au_1,v_1 \rangle + \diver{v_1} Au_1+u_2dv_2 + d\langle Av_1,u_1 \rangle+ \diver{u_1} Av_1+v_2du_2] \\ \diver(u_2v_1) + \diver(v_2u_1) \end{pmatrix}
	\\ \label{deltaexpression}
& = -\frac{1}{2}\begin{pmatrix}A^{-1}d[\langle A u_1,v_1 \rangle +\langle A v_1,u_1 \rangle -\diver{v_1}\diver{u_1}+v_2u_2 ] \\ \diver{(u_2v_1+v_2u_1)} \end{pmatrix}.
\end{align}
and
\begin{align}\label{betaexpression}
\beta & = -\frac{1}{2} \begin{pmatrix} A^{-1}[d \langle Au_1,v_1 \rangle + \diver{v_1}Au_1 +u_2dv_2 -d \langle Av_1,u_1 \rangle - \diver{u_1}Av_1- v_2du_2] \\ \diver{(u_2v_1 - v_2u_1)} \end{pmatrix}.
\end{align}
We will consider the four terms on the right-hand side of (\ref{curvatureformula}) in turn. Let $\delta_1$ and $\delta_2$ denote the two components of $\delta$. Equations (\ref{deltaexpression}) and (\ref{deltaAinvd}) yield 
\begin{align*}
\diver{\delta_1} 
    = &\; \frac{1}{2}\left(\langle A u_1,v_1 \rangle +\langle A v_1,u_1 \rangle -\diver{v_1}\diver{u_1}+v_2u_2\right)
    \\
&    - \frac{1}{2\mu(M)}\int_M \left(\langle A u_1,v_1 \rangle +\langle A v_1,u_1 \rangle -\diver{v_1}\diver{u_1}+v_2u_2\right) d\mu.
\end{align*}
Thus the first term on the right-hand side of (\ref{curvatureformula}) is given by
\begin{align}\nonumber
\llangle \delta, \delta \rrangle = &\; \frac{1}{4}\int_M (\diver{\delta_1})^2 d\mu + \frac{1}{16} \int_M (\diver(u_2v_1+v_2u_1))^2 d\mu
	\\  \label{firstterm}
= &\; \frac{1}{16} \int_M \left(\langle A u_1,v_1 \rangle +\langle A v_1,u_1 \rangle -\diver{v_1}\diver{u_1}+v_2u_2\right)^2 d\mu
 	\\ \nonumber
&  - \frac{1}{16 \mu(M)} \left\{\int_M \left(\langle A u_1,v_1 \rangle +\langle A v_1,u_1 \rangle -\diver{v_1}\diver{u_1}+v_2u_2\right) d\mu\right\}^2
	\\ \nonumber
&  + \frac{1}{16} \int_M \bigl(\diver {(u_2v_1+v_2u_1)}\bigr)^2 d\mu.
\end{align}
Using the expression (\ref{betaexpression}) for $\beta$, we deduce that the second term on the right-hand side of (\ref{curvatureformula}) is given by
\begin{align}\nonumber
\llangle [u,v], \beta \rrangle = &\; -\frac{1}{8}\int_M \Bigl\langle [u_1,v_1] , d\langle Au_1,v_1 \rangle + \diver{v_1}Au_1 +u_2dv_2 
	\\ \label{secondterm}
& \hspace{1.1cm} -d \langle Av_1,u_1 \rangle - \diver{u_1}Av_1- v_2du_2 \Bigr\rangle d\mu
	\\ \nonumber
&  - \frac{1}{8}\int_M (\langle u_1, dv_2\rangle - \langle v_1, du_2 \rangle)\diver(u_2v_1-v_2u_1) d\mu,
\end{align}
while the third term is given by
\begin{align}\label{thirdterm}
-\frac{3}{4} \llangle[u,v],[u,v] \rrangle & = -\frac{3}{16} \int_M (\diver {[u_1,v_1]})^2 d\mu -\frac{3}{16} \int_M (\langle u_1, dv_2\rangle - \langle v_1, du_2 \rangle) ^2  d\mu.
\end{align}
Moreover, 
\begin{align*}
B(u,u) & = - \begin{pmatrix} A^{-1}d\left[\langle Au_1,u_1 \rangle - \frac{1}{2}(\diver{u_1})^2 + \frac{1}{2}u_2^2\right] \\ \diver{(u_1u_2)} \end{pmatrix},
\end{align*}
so that, by (\ref{deltaAinvd}),
\begin{align*}
\diver{B_1(u,u)} = &\; \langle Au_1,u_1 \rangle - \frac{1}{2}(\diver{u_1})^2 + \frac{1}{2}u_2^2
	\\
& -\frac{1}{\mu(M)}\int_M \left(\langle Au_1,u_1^{\flat} \rangle - \frac{1}{2}(\diver{u_1})^2 + \frac{1}{2}u_2^2\right)d\mu.
\end{align*}
Thus the fourth term on the right-hand side of (\ref{curvatureformula}) is given by
\begin{align}\nonumber
 - \llangle B&(u,u),B(v,v) \rrangle = -\frac{1}{4}\int_M \diver{B_1(u,u)}\diver{B_1(v,v)}d\mu -\frac{1}{4}\int_M B_2(u,u)B_2(v,v) d\mu
	\\ \nonumber
 = & -\frac{1}{4}\int_M \left(\langle Au_1,u_1 \rangle - \frac{1}{2}(\diver{u_1})^2 + \frac{1}{2}u_2^2\right)
\left( \langle Av_1,v_1 \rangle - \frac{1}{2}(\diver{v_1})^2 + \frac{1}{2}v_2^2\right)d\mu
	\\ \label{fourthterm}
& + \frac{1}{4\mu(M)}\int_M \left(\langle Au_1,u_1 \rangle - \frac{1}{2}(\diver{u_1})^2 + \frac{1}{2}u_2^2 \right) d\mu
	\\ \nonumber
&\quad  \times \int_M \left(\langle Av_1,v_1 \rangle - \frac{1}{2}(\diver{v_1})^2 + \frac{1}{2}v_2^2 \right) d\mu 
 - \frac{1}{4}\int_M \diver{(u_1u_2)}\diver{(v_1v_2)} d\mu.
\end{align}
Adding the four contributions from (\ref{firstterm})-(\ref{fourthterm}), we find an expression for $\llangle R(u,v)v, u \rrangle$ in terms of $u_1, v_1, u_2$, and $v_2$. It is convenient to divide this expression into two terms, 
$$\llangle R(u,v)v, u \rrangle = I_1 + I_2,$$
where $I_1$ consists of all the terms that contain neither $u_2$ nor $v_2$, whereas $I_2$ consists of the remaining terms that contain $u_2$ or $v_2$. Equations (\ref{firstterm})-(\ref{fourthterm}) yield
\begin{align*}
&I_1 = \frac{1}{16} \int_M \left(\langle A u_1,v_1 \rangle +\langle A v_1,u_1 \rangle -\diver{v_1}\diver{u_1}\right)^2 d\mu
	\\
& -\frac{1}{16\mu(M)}\left\{\int_M \left(\langle A u_1,v_1 \rangle +\langle A v_1,u_1 \rangle -\diver{v_1}\diver{u_1} \right) d\mu \right\}^2
	\\
& - \frac{1}{8}\int_M \left\langle [u_1,v_1] ,d \langle Au_1,v_1 \rangle + \diver{v_1}Au_1 -d \langle Av_1,u_1 \rangle - \diver{u_1}Av_1 \right\rangle d\mu
	\\
& -\frac{3}{16} \int_M (\diver [u_1,v_1])^2 d\mu
	\\
& - \frac{1}{4}\int_M \left(\langle Au_1,u_1 \rangle- \frac{1}{2}(\diver{u_1})^2\right)\left(\langle Av_1,v_1 \rangle- \frac{1}{2}(\diver{v_1})^2 \right)d\mu
	\\
&
+\frac{1}{4\mu(M)}\left\{\int_M \left( \langle Au_1, u_1 \rangle - \frac{1}{2}(\diver{u_1})^2  \right) d\mu\right\}\left\{\int_M \left(\langle Av_1,v_1 \rangle - \frac{1}{2}(\diver{v_1})^2  \right) d\mu \right\}.
\end{align*}
Simplification using the integration by parts identity
\begin{align*}
\int_M \langle d\alpha, \beta \rangle d\mu = \int_M \langle \alpha, \delta \beta \rangle d\mu
\end{align*}
yields
\begin{align*}
I_1 = & \; \frac{1}{16} \int_M \left(\langle Au_1,v_1 \rangle +\langle A v_1,u_1 \rangle -\diver{v_1}\diver{u_1}\right)^2 d\mu
	\\
& -\frac{1}{16\mu(M)}\left(\int_M \diver{v_1}\diver{u_1} d\mu \right)^2
	\\
& +\frac{1}{8} \int_M \left(\langle Au_1,v_1 \rangle - \langle Av_1,u_1 \rangle\right)\diver{[u_1,v_1]} d\mu
	\\
& - \frac{1}{8} \int_M \left\langle [u_1,v_1], \diver{v_1}Au_1 - \diver{u_1}Av_1 \right\rangle d\mu
	\\
& -\frac{3}{16} \int_M (\diver{[u_1,v_1]})^2 d\mu
	\\
& - \frac{1}{4}\int_M \left(\langle Au_1,u_1 \rangle -\frac{1}{2}(\diver{u_1})^2 \right)\left(\langle Av_1,v_1 \rangle - \frac{1}{2}(\diver{v_1})^2\right)d\mu
	\\
& +\frac{1}{16\mu(M)}\left(\int_M (\diver{u_1})^2) d\mu\right)\left(\int_M (\diver{v_1})^2 d\mu\right).
\end{align*}
Employing the identity (\ref{lieformula}) as well as the identity
\begin{align*}
\diver{[u_1,v_1]}
& =\langle u_1, \grad {\diver{v_1}}\rangle - \langle v_1,\grad{\diver{u_1}}\rangle
=\langle Au_1, v_1 \rangle - \langle Av_1, u_1 \rangle,
\end{align*}
we can write this as
\begin{align*}
& \frac{1}{16} \int_M \left(\langle Au_1,v_1 \rangle +\langle Av_1,u_1 \rangle\right)^2 d\mu
-\frac{1}{8}\int_M \left(\langle Au_1,v_1 \rangle +\langle Av_1,u_1 \rangle\right)\diver{v_1}\diver{u_1} d\mu
	\\
&
+\frac{1}{16}\int_M (\diver{u_1})^2(\diver{v_1})^2 d\mu
-\frac{1}{16\mu(M)}\left(\int_M \diver{v_1}\diver{u_1} d\mu \right)^2
	\\
& -\frac{1}{16} \int_M \left(\langle Au_1,v_1 \rangle - \langle Av_1,u_1 \rangle\right)^2 d\mu
	\\
&
-\frac{1}{8} \int_M \left\langle u_1\diver{v_1}-v_1\diver{u_1} - (\delta(u_1^{\flat} \wedge v_1^{\flat}))^{\sharp}, 
\diver{v_1}Au_1 - \diver{u_1}Av_1 \right\rangle d\mu
	\\
&
-\frac{1}{4}\int_M \langle Au_1,u_1 \rangle \langle Av_1,v_1 \rangle d\mu
+\frac{1}{8}\int_M \langle Au_1,u_1 \rangle (\diver{v_1})^2 d\mu
	\\
&
+\frac{1}{8}\int_M \langle Av_1,v_1 \rangle (\diver{u_1})^2 d\mu
-\frac{1}{16}\int_M (\diver{u_1})^2(\diver{v_1})^2 d\mu
	\\
&+\frac{1}{16\mu(M)}\left(\int_M (\diver{u_1})^2 d\mu\right)\left(\int_M (\diver{v_1})^2 d\mu\right).
\end{align*}
Finally, using that
\begin{align*}
& \int_M \left\langle \delta(u_1^{\flat} \wedge v_1^{\flat}), 
\diver{v_1}Au_1 - \diver{u_1}Av_1 \right\rangle d\mu 
= 2 \int_M \left\langle u_1^\flat \wedge v_1^\flat, Au_1 \wedge Av_1 \right \rangle d\mu
	\\
&\hspace{4.2cm} = 2\int_M \left[\langle u_1, Au_1\rangle\langle v_1, Av_1 \rangle - \langle u_1, Av_1 \rangle\langle v_1, Au_1 \rangle\right] d\mu,
\end{align*}
we arrive after several cancellations at the following expression for $I_1$:
\begin{align}\label{I1final}
I_1 =\frac{1}{16\mu(M)}\left(\int_M (\diver{u_1})^2 d\mu\right)\left(\int_M (\diver{v_1})^2 d\mu\right)
- \frac{1}{16\mu(M)}\left(\int_M \diver{v_1}\diver{u_1} d\mu \right)^2.
\end{align}

On the other hand, the terms in equations (\ref{firstterm})-(\ref{fourthterm}) that contain $u_2$ or $v_2$ yield
\begin{align} \nonumber
I_2 = & \;\frac{1}{16} \int_M \left( 2(\langle Au_1,v_1 \rangle + \langle Av_1, u_1 \rangle 
- \diver{v_1}\diver{u_1})u_2v_2 + (u_2v_2)^2 \right) d\mu
	\\ \nonumber
& -\frac{1}{8\mu(M)}\left(\int_M (\langle Au_1,v_1\rangle +\langle Av_1,u_1 \rangle -\diver{v_1}\diver{u_1}) d\mu\right)\left(\int_M u_2v_2 d\mu\right)
	\\ \nonumber
& -\frac{1}{16\mu(M)} \left(\int_M u_2v_2d\mu\right)^2
 +\frac{1}{16} \int_M (\diver{(u_2v_1+v_2u_1)})^2 d\mu
	\\ \nonumber
& -\frac{1}{8}\int_M \left\langle [u_1,v_1], u_2dv_2-v_2du_2 \right\rangle d\mu
	\\ \nonumber
& -\frac{1}{8}\int_M (\langle u_1, dv_2 \rangle - \langle v_1, du_2 \rangle)\diver(u_2v_1-v_2u_1) d\mu
	\\ \label{longeq}
& -\frac{3}{16}\int_M \left(\langle u_1, dv_2 \rangle - \langle v_1, du_2 \rangle\right)^2 d\mu
	\\ \nonumber
&- \frac{1}{8}\int_M \left(\langle Au_1, u_1 \rangle - \frac{1}{2}(\diver{u_1})^2\right)v_2^2 d\mu
	\\ \nonumber
&
- \frac{1}{8}\int_M \left(\langle Av_1,v_1 \rangle - \frac{1}{2}(\diver{v_1})^2\right)u_2^2 d\mu
-\frac{1}{16}\int_M (u_2v_2)^2 d\mu
	\\ \nonumber
&+ \frac{1}{8\mu(M)}\left\{\int_M \left(\langle Au_1,u_1 \rangle - \frac{1}{2}(\diver{u_1})^2\right) d\mu\right\}\left(\int_M v_2^2 d\mu\right)
	\\ \nonumber
&
+ \frac{1}{8\mu(M)}\left\{\int_M \left(\langle Av_1, v_1 \rangle - \frac{1}{2}(\diver{v_1})^2\right) d\mu\right\}\left(\int_M u_2^2 d\mu\right)
	\\ \nonumber
&+\frac{1}{16\mu(M)}\left(\int_M u_2^2 d\mu\right)\left(\int_M v_2^2 d\mu\right)
-\frac{1}{4}\int_M \diver(u_1u_2) \diver(v_1v_2) d\mu.
\end{align}
We use the identity
$$\diver(u_2v_1+v_2u_1) = u_2\diver{v_1}+\langle du_2,v_1 \rangle + v_2\diver{u_1}+\langle dv_2, u_1 \rangle$$
in the third line, and the identity (\ref{lieformula}) in the fourth line of (\ref{longeq}). Moreover, we use the identity
$$\diver(u_2v_1-v_2u_1) = u_2\diver{v_1} + \langle du_2,v_1 \rangle - v_2\diver{u_1} - \langle dv_2, u_1 \rangle$$
in the fifth line of (\ref{longeq}) and combine the result with the sixth line. After simplification this leads to
\begin{align*}
I_2 =&\; \frac{1}{8} \int_M u_2v_2\langle Au_1,v_1 \rangle d\mu
+\frac{1}{8}\int_M u_2v_2\langle Av_1,u_1 \rangle d\mu
	\\
& -\frac{1}{8}\int_M u_2v_2\diver{v_1}\diver{u_1} d\mu
+\frac{1}{16}\int_M (u_2v_2)^2 d\mu
	\\
&
-\frac{1}{8\mu(M)}\left(\int_M \diver{v_1}\diver{u_1} d\mu\right)\left(\int_M u_2v_2 d\mu\right)
-\frac{1}{16\mu(M)}\left(\int_M u_2v_2 d\mu \right)^2
	\\
&
+\frac{1}{16} \int_M \biggl\{ u_2^2(\diver{v_1})^2+ \langle du_2,v_1 \rangle ^2 + v_2^2(\diver{u_1})^2+ \langle dv_2, u_1 \rangle^2
	\\
&\hspace{1.8cm}+ 2u_2\diver{v_1}\langle du_2, v_1 \rangle+ 2u_2\diver{v_1}v_2\diver{u_1} + 2u_2\diver{v_1}\langle dv_2, u_1 \rangle
		\\
&\hspace{1.8cm}+ 2\langle du_2,v_1 \rangle  v_2\diver{u_1}+ 2\langle du_2,v_1 \rangle \langle dv_2,u_1 \rangle +2 v_2\diver{u_1} \langle dv_2,u_1 \rangle   \biggr\} d\mu
	\\
&
- \frac{1}{8}\int_M \Bigl\{\langle u_1\diver{v_1}, u_2dv_2\rangle - \langle u_1\diver{v_1}, v_2du_2\rangle
-\langle v_1\diver{u_1}, u_2dv_2\rangle 
	\\
& \hspace{1.7cm} + \langle v_1\diver{u_1}, v_2du_2 \rangle 
 - \left\langle \delta(u_1^{\flat}\wedge v_1^{\flat}), u_2dv_2-v_2du_2 \right\rangle   \Bigr\} d\mu
	\\
&
-\frac{1}{8}\int_M \Bigl\{ u_2\diver{v_1}\langle u_1,  dv_2\rangle 
- v_2\diver{u_1}\langle u_1, dv_2 \rangle
	\\ 
& \hspace{1.7cm} - u_2\diver{v_1}\langle v_1,  du_2 \rangle
+ v_2\diver{u_1} \langle v_1, du_2 \rangle \Bigr\} d\mu
	\\
&
-\frac{1}{16}\int_M \left(\langle u_1, dv_2\rangle ^2+ \langle v_1, du_2 \rangle ^2
 -2\langle u_1, dv_2\rangle \langle v_1, du_2 \rangle \right) d\mu
	\\
&
-\frac{1}{8}\int_M v_2^2\langle Au_1,u_1 \rangle d\mu
+\frac{1}{16}\int_M v_2^2(\diver{u_1})^2 d\mu
	\\
& -\frac{1}{8}\int_M u_2^2\langle Av_1,v_1 \rangle d\mu
 +\frac{1}{16}\int_M u_2^2(\diver{v_1})^2 d\mu
-\frac{1}{16}\int_M (u_2v_2)^2 d\mu
	\\
&
+\frac{1}{16\mu(M)}\left(\int_M(\diver{u_1})^2 d\mu\right)\left(\int_M v_2^2 d\mu\right)
+\frac{1}{16\mu(M)}\left(\int_M(\diver{v_1})^2 d\mu\right)\left(\int_M u_2^2 d\mu\right)
	\\
&
+\frac{1}{16\mu(M)}\left(\int_M u_2^2 d\mu\right)\left(\int_M v_2^2 d\mu\right)
	\\
&
-\frac{1}{4}\int_M \Bigl\{(u_2v_2\diver{u_1}\diver{v_1}+u_2\diver{u_1}\langle v_1, dv_2\rangle
	\\
& \hspace{1.7cm}+v_2\diver{v_1}\langle u_1, du_2\rangle
+\langle u_1,du_2\rangle \langle v_1, dv_2\rangle \Bigr\} d\mu.
\end{align*}
Using the identity
\begin{align*}
 \int_M \left\langle \delta(u_1^{\flat} \wedge v_1^{\flat}), u_2dv_2 - v_2du_2 \right\rangle &d\mu 
 = 2 \int_M \left\langle u_1^\flat \wedge v_1^\flat, du_2 \wedge dv_2 \right \rangle d\mu
	\\
& = 2\int_M \left[\langle u_1, du_2\rangle\langle v_1, dv_2 \rangle - \langle u_1, dv_2 \rangle\langle v_1, du_2 \rangle\right]d\mu
\end{align*}
as well as integration by parts, we arrive at the following expression for $I_2$:
\begin{align}\nonumber
I_2 = &-\frac{1}{8\mu(M)}\left(\int_M \diver{v_1}\diver{u_1} d\mu\right)\left(\int_M u_2v_2 d\mu\right)
-\frac{1}{16\mu(M)}\left(\int_M u_2v_2d\mu\right)^2
	\\ \nonumber
& +\frac{1}{16\mu(M)}\left(\int_M(\diver{u_1})^2 d\mu\right)\left(\int_M v_2^2 d\mu\right)
	\\ \label{I2final}
& +\frac{1}{16\mu(M)}\left(\int_M(\diver{v_1})^2 d\mu\right)\left(\int_M u_2^2 d\mu\right)
    \\ \nonumber
& +\frac{1}{16\mu(M)}\left(\int_M u_2^2 d\mu\right)\left(\int_M v_2^2 d\mu\right).
\end{align}

Finally, addition of the expressions (\ref{I1final}) and (\ref{I2final}) for $I_1$Ê and $I_2$ yields
\begin{align*}
 \llangle R(u,v)v,u \rrangle 
=&\; \frac{1}{16\mu(M)} \left(\int_M(\diver{u_1})^2 d\mu+ \int_M u_2^2d\mu\right)\left(\int_M(\diver{v_1})^2 d\mu+ \int_M v_2^2d\mu\right)
    \\
& -\frac{1}{16\mu(M)}\left(\int_M \diver{u_1}\diver{v_1}d\mu+ \int_M u_2v_2\right)^2
    \\
= &\;  \frac{\llangle u, u \rrangle \llangle v, v \rrangle - \llangle u, v\rrangle^2}{\mu(M)}.
\end{align*}
This proves (\ref{Rexpression}) and hence completes the proof of Theorem \ref{mainth}.

\appendix
\section{The curvature of a homogeneous space}\label{curvatureapp}\nequation
\renewcommand{\theequation}{A.\arabic{equation}}
In this appendix we derive the curvature formula (\ref{curvatureformula}) for a homogeneous space with a right-invariant metric. 

Suppose that $G$ is a Lie group and that $H \subset G$ is a (closed) subgroup of $G$. Let $M = G/H$ denote the homogeneous space of right cosets. The Lie algebras of $G$ and $H$ are denoted by $\mathfrak{g}$ and $\mathfrak{h}$, and we let $\mathfrak{m} \subset \mathfrak{g}$ be a linear subspace such that $\mathfrak{g} = \mathfrak{h} \oplus \mathfrak{m}$. Letting $\pi:G \to G/H$ denote the quotient map, we find that $\pi_*$ maps $\mathfrak{m}$ isomorphically onto $T_{\pi(e)}M$, where $e$ denotes the identity element in $G$.
Suppose that $\langle \cdot, \cdot \rangle$ is a right-invariant (degenerate) metric on $G$ whose restriction to $\mathfrak{m}$ is positive definite and whose kernel is $\mathfrak{h}$, i.e. $\langle u, v \rangle = 0$ for all $v \in \mathfrak{g}$ iff $u \in \mathfrak{h}$. We assume that the restriction of the metric $\langle \cdot, \cdot \rangle$ to $\mathfrak{m}$ is $\Ad(H)$-invariant so that $\langle \cdot, \cdot \rangle$ descends to a right-invariant metric on $G/H$ which will be denoted by $\met( \cdot, \cdot)$.  

\begin{remark}\upshape
1. We do {\it not} require that $\mathfrak{m}$ be an $\Ad(H)$-invariant subspace. In particular, the homogeneous space $G/H$ does not have to be reductive.

2. For the homogeneous space of Theorem \ref{mainth}, we have 
$$\mathfrak{g} = \mathfrak{X}(M) \times C^\infty(M), \qquad
\mathfrak{h} = \{(v, 0) \in \mathfrak{g} \, | \, \diver{v} = 0\},$$ 
and we may choose $\mathfrak{m} = \{((df)^\sharp, g) \,|\, f, g \in C^\infty(M)\}$. 
\end{remark}

For any $X \in \mathfrak{g}$, we define the vector field $\tilde{X}\in \mathfrak{X}(G/H)$ as the push-forward by $\pi$ of $X^L$, where $X^L$ denotes the unique left-invariant vector field on $G$ whose value at $e$ is $X$. Thus,
$$\tilde{X}_{\pi(g)} = \frac{d}{dt}\bigg|_{t=0} \pi(g e^{tX}) = \pi_*X^L_g, \qquad g \in G.$$
where $e^{tX}$ is the one-parameter subgroup in $G$ generated by $X$. The flow $\Phi_X^t:M \to M$ of $\tilde{X}$ is given by
$$\Phi_X^t(\pi(g)) = \pi(g e^{tX}), \qquad g \in G, \; t \in \R.$$
If $U \in TG$, then $(\Phi_X^t)_*\pi_*U = \pi_* R_{e^{tX}*} U$, where $R_g$ denotes right multiplication by $g \in G$. Hence, since $\langle \cdot, \cdot \rangle$ is right invariant,
$$\met((\Phi_X^t)_*\pi_* U, (\Phi_X^t)_*\pi_* V)
= \langle R_{e^{tX}*} U, R_{e^{tX}*} V \rangle
= \langle U, V \rangle = \met(\pi_*U, \pi_*V),$$
whenever $U, V \in T_gG$, that is, the flow of $\tilde{X}$ consists of isometries. This implies that $\tilde{X}$ is a Killing field. In particular, $\tilde{X}$ satisfies
\begin{align}\label{killing1}
\tilde{X} \met(Y, Z) = \met([\tilde{X}, Y], Z) + \met(Y, [\tilde{X}, Z])
\end{align}
and
\begin{align}\label{killing2}
  \met(\nabla_Y \tilde{X}, Z) + \met(\nabla_Z \tilde{X}, Y) = 0
\end{align}  
for all vector fields $Y, Z$ on $G/H$. 

Note that
\begin{align}\label{liestar}
[\tilde{X}, \tilde{Y}] = -\widetilde{[X, Y]}, \qquad X, Y \in \mathfrak{g},
\end{align}
where the bracket on the left-hand side is the Lie bracket of vector fields on $M$ and the bracket on the right-hand side is the Lie bracket in $\mathfrak{g}$ induced by {\it right}-invariant vector fields (this bracket is minus the bracket induced by left-invariant vector fields). 
Indeed, 
$$[\tilde{X}, \tilde{Y}] = [\pi_* X^L, \pi_* Y^L] = \pi_*[X^L, Y^L] = -\pi_*[X, Y]^L = -\widetilde{[X, Y]}.$$

Let $\nabla$ denote the Levi-Civita connection associated with the Riemannian metric $g$ on $M = G/H$. Then, for any vector fields $X,Y,Z$ on $M$,
\begin{align*}
2\met(\nabla_X Y, Z) = &\; X\met(Y, Z) + Y\met(Z, X) - Z\met(X, Y) 
	\\
&+ \met(Z, [X,Y]) - \met(Y, [X, Z]) - \met(X, [Y, Z]).
\end{align*}
In the particular case when $X,Y,Z$ are Killing vector fields, (\ref{killing1}) implies that
\begin{align}\label{nablakilling}
  2\met(\nabla_X Y, Z) = \met([X, Y], Z) + \met([X, Z], Y) + \met(X, [Y, Z]).
\end{align}

\begin{lemma}
Let $X,Y \in \mathfrak{m}$. Then
\begin{align}\label{nablaXYstar}
(\nabla_{\tilde{X}}\tilde{Y})_{\pi(e)} = -\frac{1}{2}\pi_*([X, Y] + B(X,Y) + B(Y,X)),
\end{align}
where $B(X,Y) \in \mathfrak{g}$ is defined by
\begin{align}\label{BXYZdef}
\langle B(X, Y), Z \rangle = \langle X, [Y, Z] \rangle.
\end{align}
\end{lemma}
\begin{remark}\upshape
Equation (\ref{BXYZdef}) only determines $B(X,Y)$ up to addition by an element in $\mathfrak{h}$. We can fix this freedom by requiring that $B(X,Y) \in \mathfrak{m}$. Alternatively, we can just note that $\pi_*B(X,Y)$ is uniquely determined and so the right-hand side of (\ref{nablaXYstar}) is well-defined. 
\end{remark}
\proofbegin
Let $X, Y, Z \in \mathfrak{m}$. Then (\ref{nablakilling}) and (\ref{liestar}) give
$$2\met(\nabla_{\tilde{X}} \tilde{Y}, \tilde{Z}) = -\met(\widetilde{[X, Y]}, \tilde{Z}) - \met(\widetilde{[X, Z]}, \tilde{Y}) - \met(\tilde{X}, \widetilde{[Y, Z]}).$$
Evaluating this equation at the point $\pi(e)$, we find
$$2\met\Bigl((\nabla_{\tilde{X}} \tilde{Y})_{\pi(e)} + \frac{1}{2}\widetilde{[X, Y]}_{\pi(e)}, \tilde{Z}_{\pi(e)}\Bigr) = -\langle [X, Z], Y \rangle - \langle X, [Y, Z] \rangle,$$
and this equality gives (\ref{nablaXYstar}). 
\proofend

\begin{proposition}
Let $X, Y \in \mathfrak{m} \simeq T_{\pi(e)}M$. The curvature tensor $R$ of $M=G/H$ satisfies
\begin{align}\label{RXYYX}
\met( R(X,Y)Y, X)
 = \langle \delta, \delta \rangle + \langle [X,Y], \beta \rangle - \frac{3}{4}\langle [X,Y], [X,Y]\rangle - \langle B(X,X), B(Y,Y) \rangle,
\end{align}
where $\delta$ and $\beta$ are defined by
$$\delta = \frac{1}{2}(B(X,Y) + B(Y,X)), \qquad \beta = \frac{1}{2}(B(X,Y) - B(Y,X)).$$
\end{proposition}
\proofbegin
We compute
$$\met(R(\tilde{X},\tilde{Y})\tilde{X}, \tilde{Y}) = \met(\nabla_{\tilde{X}}\nabla_{\tilde{Y}}\tilde{X}, \tilde{Y}) - \met(\nabla_{\tilde{Y}}\nabla_{\tilde{X}}\tilde{X}, \tilde{Y}) - \met(\nabla_{[\tilde{X},\tilde{Y}]}\tilde{X}, \tilde{Y}).$$
This gives, using (\ref{killing2}) in the last term,
\begin{align*}
\met(R(\tilde{X},\tilde{Y})\tilde{X}, \tilde{Y}) = &\; \tilde{X} \met(\nabla_{\tilde{Y}}\tilde{X}, \tilde{Y})  - \met(\nabla_{\tilde{Y}}\tilde{X}, \nabla_{\tilde{X}} \tilde{Y}) 
	\\
& - \tilde{Y} \met(\nabla_{\tilde{X}}\tilde{X}, \tilde{Y}) + \met(\nabla_{\tilde{X}} \tilde{X}, \nabla_{\tilde{Y}} \tilde{Y}) + \met(\nabla_{\tilde{Y}}\tilde{X}, [\tilde{X},\tilde{Y}]).
\end{align*}
Using the following two equations, which are consequences of (\ref{nablakilling}),
\begin{align*}
& 2\met(\nabla_{\tilde{Y}}\tilde{X}, \tilde{Y}) = \met([\tilde{Y}, \tilde{X}], \tilde{Y}) + \met([\tilde{Y}, \tilde{Y}], \tilde{X}) + \met(\tilde{Y}, [\tilde{X}, \tilde{Y}]) = 0,
	\\
& 2\met(\nabla_{\tilde{X}}\tilde{X}, \tilde{Y}) =  \met([\tilde{X}, \tilde{X}], \tilde{Y}) + \met([\tilde{X}, \tilde{Y}], \tilde{X}) + \met(\tilde{X}, [\tilde{X}, \tilde{Y}]) = 2\met([\tilde{X}, \tilde{Y}], \tilde{X}),
\end{align*}
as well as the fact that $\nabla$ is torsion-free, i.e. $\nabla_{\tilde{X}}\tilde{Y} - \nabla_{\tilde{Y}}\tilde{X} = [\tilde{X}, \tilde{Y}]$, we find
$$\met(R(\tilde{X},\tilde{Y})\tilde{X}, \tilde{Y}) = - \tilde{Y}\met([\tilde{X}, \tilde{Y}], \tilde{X})+ \met(\nabla_{\tilde{X}} \tilde{X}, \nabla_{\tilde{Y}} \tilde{Y}) - \met(\nabla_{\tilde{Y}}\tilde{X}, \nabla_{\tilde{Y}}\tilde{X}).$$
Equation (\ref{killing1}) implies that 
$$\tilde{Y}\met([\tilde{X}, \tilde{Y}], \tilde{X}) = \met([\tilde{Y}, [\tilde{X}, \tilde{Y}]], \tilde{X}) + \met([\tilde{X}, \tilde{Y}], [\tilde{Y}, \tilde{X}]),$$
so that, in view of (\ref{nablaXYstar}),
\begin{align*}
\met(R(\tilde{X},\tilde{Y})\tilde{X}, \tilde{Y}) = & - \met([\tilde{Y}, [\tilde{X}, \tilde{Y}]], \tilde{X}) - \met([\tilde{X}, \tilde{Y}], [\tilde{Y}, \tilde{X}])  + \langle B(X,X), B(Y,Y) \rangle
	\\
& - \frac{1}{4}\langle[Y, X], [Y, X]Ê\rangle
- \langle [Y, X], \delta \rangle - \langle \delta, \delta \rangle.
\end{align*}
Finally, the relation (\ref{liestar}) yields
\begin{align*}
\met(R(\tilde{X},\tilde{Y})\tilde{X}, \tilde{Y}) = & 	- \langle B(X, Y), [X, Y] \rangle + \langle [X, Y], [X, Y] \rangle + \langle B(X,X), B(Y,Y) \rangle 
	\\
& - \frac{1}{4}\langle[X, Y], [X, Y]Ê\rangle
+ \langle [X, Y], \delta \rangle - \langle \delta, \delta \rangle.
\end{align*}
This is (\ref{RXYYX}).
\proofend

\bigskip
\noindent
{\bf Acknowledgement} {\it J.L. acknowledges support from the EPSRC, UK.}\\

\bibliographystyle{plain}
\bibliography{is}

\begin{thebibliography}{99}
\small

\bibitem{A1966}
V.I. Arnold,
S\"ur la g\'eometrie diff\'erentielle des groupes de Lie de dimension infinie et
ses application \`a l'hydrodynamique des fluides parfaits,
{\it Ann. Inst. Fourier (Grenoble)} {\bf 16} (1966), 319--361.

\bibitem{Abook}
V.I. Arnold, {\it Mathematical methods of classical mechanics}, 2nd edition, Graduate Texts in Mathematics 60, Springer-Verlag, New York, 1989.

\bibitem{CI2008}
A. Constantin and R. Ivanov, On an integrable two-component Camassa-Holm shallow water system, {\it Phys. Lett. A} {\bf 372} (2008), 7129--7132.

\bibitem{EM1970}
D. Ebin and J.E. Marsden,
Groups of diffeomorphisms and the motion of an incompressible fluid,
{\it Ann. Math.} {\bf 92} (1970), 102--163.

\bibitem{EKL2011}
J. Escher, M. Kohlmann, and J. Lenells, The geometry of the two-component Camassa-Holm and Degasperis-Procesi equations, {\it J. Geom. Phys.} {\bf 61} (2011), 436--452.

\bibitem{HS1991}
J.K. Hunter and R. Saxton, Dynamics of director fields,
{\it SIAM J. Appl. Math.} {\bf 51} (1991), 1498--1521.

\bibitem{KMLP2011}
B. Khesin, J. Lenells, G. Misio\l ek, and S.C. Preston, Geometry of diffeomorphism groups, complete integrability and geometric statistics, arXiv:1105.0643.

\bibitem{KM2003}
B. Khesin and G. Misio\l ek, Euler equations on homogeneous spaces and Virasoro orbits, {\it Adv. Math.} {\bf 176} (2003), 116--144.

\bibitem{KN1969}
S. Kobayashi and K. Nomizu, {\it Foundations of Differential Geometry II}, John Wiley \& Sons, New York, 1969.

\bibitem{K2011}
M. Kohlmann, The curvature of semidirect product groups associated with 2HS and 2$\mu$HS, preprint.

\bibitem{Lang}
S. Lang, {\it Differential and Riemannian Manifolds}, 3rd ed., Springer-Verlag, New York, 1995.

\bibitem{Lsphere}
J. Lenells,
The Hunter-Saxton equation describes the geodesic flow on a sphere,
{\it J. Geom. Phys.} {\bf 57} (2007), 2049--2064.

\bibitem{M1998}
G. Misio\l ek, A shallow water equation as a geodesic flow on the
        Bott-Virasoro group, {\it J. Geom. Phys.} {\bf 24} (1998),
        203--208.
        
\bibitem{R1997}
S. Rosenberg, {\it The Laplacian on a Riemannian manifold}, London Mathematical Society Student Texts {\bf 31}. Cambridge University Press, Cambridge, 1997.

\bibitem{V2008}
C. Vizman, Geodesic equations on diffeomorphism groups, {\it SIGMA Symmetry Integrability Geom. Methods Appl.} {\bf 4} (2008), Paper 030, 22 pp.

\bibitem{W2009}
M. Wunsch, On the Hunter-Saxton system, {\it Disc. Cont. Dyn. Syst. Ser. B} {\bf 12} (2009), 647--656. 

\end{thebibliography}

\end{document}